\documentclass[11pt,titlepage]{article}
 \usepackage{epic}
 \usepackage{curves}
 \usepackage{epsf}
 \usepackage{amssymb}
 \usepackage{latexsym}

 \date{2nd of October 2000}

 \newcommand{\grad}{\mathop{\rm grad}\nolimits}
 
 \newcommand{\im}{\mathop{\rm im}\nolimits}
 \newcommand{\R}{{\mathbb R}}
 \newcommand{\Z}{{\mathbb Z}}
 \newcommand{\RP}{{\mathbb R}{\rm P}}

 \newtheorem{theorem}{Theorem}
 \newtheorem{lemma}{Lemma}[section]
 \newtheorem{note}{Remark}
 
 \newtheorem{corollary}{Corollary}
 
 \newcommand{\beginproof}{\begin{trivlist}
 \rm\item[\hspace{\labelsep}{ P\,r\,o\,o\,f. }]}
 \newcommand{\proofend}{$\Box$\end{trivlist}}
 \makeatletter
 \renewcommand{\@begintheorem}[2]{\begin{trivlist}\it\item[\hspace{\labelsep}
  {\sc #1\ #2}]}
 \renewcommand{\@endtheorem}{\end{trivlist}}
 \makeatother

\begin{document}
 \begin{titlepage}
 \begin{center}

 {\Large\sc Moscow State University

 \bigskip
 Department of Mechanics and Mathematics

 }
 \vspace{5cm}

 {\LARGE\sc
 On the lower bounds for the number of\\[1mm]
 periodic billiard trajectories in manifolds \\[3.7mm]
 embedded in Euclidean space
 }

 \vspace{3cm}

 \LARGE\rm 
 
 Fedor Duzhin
 
\vspace{.7cm}

 \Large
M.S. thesis

 \vspace{3cm}
 \end{center}

 \begin{flushright}

 Scientific advisor\\professor S.~Gusein-Zade

 \end{flushright}

 \vfill

 \begin{center}
 {\LARGE\sc Moscow 2000}
 \end{center}
 \end{titlepage}

 \section{Introduction}

 We shall study the lower bounds for the number of the periodic
 billiard trajectories in manifolds embedded in Euclidean space.
 A $p$-periodic billiard trajectory is a closed polygon consisting of
 $p$ segments all of whose vertices belong to the given manifold and,
 at every vertex, the two angles formed by the line and the manifold
 are equal (the exact definition will be given later).
 The first who considered this problem was George Birkhoff.
 He proved the following fact in \cite{birkhoff}.
 Suppose $p$ is an odd prime,
 $M$ is a strictly convex smooth closed curve.
 Then there exist in $M$ at least two periodic
 billiard trajectories for each rotation number
 from $1$ to $\frac{p-1}{2}$.
 Ivan Babenko studied the billiards in a $2$-dimensional sphere,
 but his paper \cite{babenko} contains an error.
 Peter Pushkar in \cite{pushkar} solved the problem for $p=2$.
 He showed that in generic case
 in an $m$-dimensional manifold there are at least
 $\frac{B^2+(m-1)B}{2}$ diameters, i. e., $2$-periodic
 trajectories, where $B$ is the sum of Betty numbers
 modulo $2$ of the given manifold.
 Finally,
 Michael Farber and Serge Tabachnikov proved
 in \cite{farber} that for $m\ge 3$ in an $m$-dimensional
 sphere there exist at least
 $\log_2(p-1)+m$ ($m(p-1)$ in a generic case)
 $p$-periodic billiard trajectories.

 In Section 2 we show how one can apply
 Morse theory to study
 periodic billiard trajectories.
 Section 3 contains the generalized Birkhoff theorem.
 In Section 4 we prove the Farber-Tabachnikov
 estimate for generic small perturbations
 of the standard round $m$-sphere for any $m$.
 In Section 5 we study
 $3$-periodic billiard trajectories in a $2$-dimensional sphere.
 At last in Section 6 we find the rough estimate
 for the number of $3$-periodic billiard trajectories
 in any manifold.

 \section{Morse theory of periodic billiard trajectories}

 Suppose $M$ is a smooth closed $m$-dimensional manifold
 embedded in Euclidean space $\R ^n$,
 $p$ is an odd prime. An ordered set of points
 $(x_1,\dots ,x_p)\in M\times\dots\times M$ is called
 a $p$-periodic billiard trajectory if
 $$
 \frac{x_i-x_{i+1}}{\|x_i-x_{i+1}\|}+
 \frac{x_i-x_{i-1}}{\|x_i-x_{i-1}\|}\ \bot\ T_{x_i}M
 $$
 for any $x_i$, $i\in \Z _p$.
 We consider sets of points
 $(x_1,\dots ,x_p)\in M\times\dots\times M$
 up to the action of the dihedral group
 $D_p$ in the $p$th Cartesian power of $M$.
 This action is generated by
 the cyclic permutation and the reflection:
 $$
 (x_1, x_2, \ldots, x_p) \rightarrow (x_2, \ldots, x_p, x_1)
 $$
 $$
 (x_1, x_2, \ldots, x_p) \rightarrow (x_p, x_{p-1}, \ldots, x_1)
 $$

 Suppose $\Delta =\{(x,x,x_3,\dots,x_p)\}\subset M\times\dots\times M/D_p$.
 Let $f$ be the length function of a closed polygon:
 $$
 f(x_1,\dots,x_p)=\sum_{i\in\Z_p}\rho(x_i,x_{i+1}).
 $$
 Thus $\Delta$ is the set of all points where $f$ is not smooth.
 It is clear that
 $p$-periodic billiard trajectories are
 critical points of $f$.
 Indeed,
 $$
 \frac{\partial f}{\partial x_i}=
 \frac{x_i-x_{i+1}}{\|x_i-x_{i+1}\|}+
 \frac{x_i-x_{i-1}}{\|x_i-x_{i-1}\|},\ x_i\in {\R}^n.
 $$
 The derivative along any tangent vector vanishes
 if and only if the gradient is orthogonal to
 the tangent space.

 In this paper we consider the general case:
 $f$ is a Morse function outside of $\Delta$.
 By $BT_p(M)$ denote
 the minimal number of $p$-periodic
 billiard trajectories in $M$.

 \medskip

 \begin{lemma} Let $M$ be
 a closed Riemannian manifold, $\dim M=m$.
 Then there exists
 $\varepsilon >0$ such that
 the following conditions hold for any $x\in M$:

 1) The solid sphere $B_{\varepsilon}(x) \subset M$
 is diffeomorphic to the disk
 $D_{\varepsilon}(0) \subset \R ^m$;
 this diffeomorphism maps geodesics passing through $x\in M$
 to straight lines passing through $0\in\R ^m$;
 angles between geodesics preserve.

 2) Suppose
 $y_1(t), y_2(t),\dots, y_k(t)$
 are mutually different geodesics such that
 $y_1(0)=y_2(0)=\dots=y_k(0)=x$.
 Let $\alpha _1,\dots ,\alpha _k>0$ be real numbers. Put
 $$
 h(t)=\rho(x,y_1(\alpha_1t))+
 \sum_{i=1}^{k-1}\rho(y_i(\alpha_it),y_{i+1}(\alpha_{i+1}t)),\ t>0
 $$
 If $\alpha _1t,\dots ,\alpha _kt<\varepsilon $,
 then $h'(t)\ge\alpha_k$.
 \end{lemma}

 \beginproof Since $M$ is compact, it is sufficient to find the required
 $\varepsilon$ only for one point $x\in M$.
 We introduce coordinates in a neighborhood of the point $x$
 along geodesics passing through $x$.
 In these coordinates
 $$
 \rho (y,z)=\|y-z\|+O(\|y-z\|^2).
 $$
 Thus for $h(t)$ we have
 $$
  h'(t)=\alpha_1+
 \sum_{i=1}^{k-1}\|\alpha_i{\bf a}_i-\alpha_{i+1}{\bf a}_{i+1}\|
 +O(t)\ge\alpha_k
 $$
 by the triangle inequality applied to the points
 $0,\alpha_1{\bf a}_1,\ldots,\alpha_k{\bf a}_k$,
 where ${\bf a}_i$
 are the directing vectors of the straight lines
 corresponding to the geodesics $y_i(t)$.
 Thus we find the required neighborhood of $x$.
 \proofend

 \medskip

 \begin{note} We can choose
 $\varepsilon$ so that the 2nd condition holds for each
 $k=1,\ldots,K$.
 \end{note}

 \medskip

 \begin{theorem}
 The minimal number of p-periodic billiard trajectories satisfies
 $$
 BT_p(M)\ge \sum_{q=1}^{mp}\dim H_q(M\times\dots\times M/D_p,\Delta;\Z_2).
 $$
 \end{theorem}

 \beginproof We construct a function $g$
 on $X=M\times\dots\times M/D_p$ such that the following
 conditions hold:

 1. $g\ge 0$,

 2. $g$ is smooth outside of $\Delta$,

 3. Critical points of $g$ are the same as those of $f$,

 4. $\Delta=\{g=0\}$.
 \\Then our theorem follows from Morse theory.
 Indeed, a small neighborhood of $\Delta$
 in $X/\Delta$ is contractible, that's why
 we can construct a cell space using the function $g$
 in the same way as using any Morse function.
 Further, we apply Morse inequalities and reduce
 the relative homology to the absolute:
 $$
 \tilde H_*(X/\Delta)=H_*(X,\Delta).
 $$

 Suppose $\varphi(t)$ is a smooth function such that
 ${0\le\varphi\le 1}$, ${\varphi\mid _{(-\infty,0]}\equiv 0}$,
 ${\varphi\mid _{[\varepsilon,+\infty)}\equiv 1}$,
 ${\varphi'\mid _{(0,\varepsilon)}>0}$.
 Let us show that the function
 $$
 g(x_1,\dots,x_p)=
 f(x_1,\dots,x_p)
 \left(\prod_{i\in\Z_p}\varphi(\rho(x_i,x_{i+1}))\right)
 $$
 is required if $\varepsilon$ is small enough.

 By definition, put
 $\Delta^{(p-k)}=\{(x_1,\dots,x_1,x_{k+1},\dots,x_p\}\subset \Delta$.
 Then $\Delta^{(0)}\subset\ldots\subset\Delta^{(p-2)}=\Delta$.
 First we find
 $\varepsilon$ for a neighborhood of $\Delta^{(0)}$.
 Suppose $\varepsilon$ is chosen by the previous remark
 with $K=p-1$.
 Then $f'(t)>0$ for $t<\varepsilon$.
 If necessary we decrease $\varepsilon$ so that
 all critical points of $f$ are outside of
 the $\varepsilon$-neighborhood of $\Delta$. Then
 $$
 U_{\varepsilon}^0(\Delta^{(0)})=\{(x_1,\dots,x_p):
 \rho(x_1,x_2),\ldots,\rho(x_1,x_p)<\varepsilon\}
 $$
 is a small neighborhood of $\Delta^{(0)}$.

 We can obtain any point $(x_1,x_2,\dots,x_p)$ in this neighborhood
 if we fix
 $x_1$, emit $p-1$ geodesics from $x_1$,
 and put points $x_2,\dots,x_p$ on these geodesics.

 Let $x_j(\alpha_jt), j\ge 2$ be geodesics
 passing through $x_1$, $t$ is a natural parameter, i.~e.,
 $\|\frac{d}{dt}x_j(\alpha_jt)\|=\alpha_j$. Then $g(x_1,\dots,x_p)=g(t)$.
 Compute the derivative:
 $$
 g'(t)=\left(\prod_{i\in\Z_p}\varphi(\rho(x_i,x_{i+1}))\right)'f(t)+
 \left(\prod_{i\in\Z_p}\varphi(\rho(x_i,x_{i+1}))\right)f'(t).
 $$
 Inequalities
 $$
 \varphi>0,\ f>0,\ \rho(x_i,x_{i+1})>0,
 $$
 $$
 \varphi'\ge 0,\ f'(t)>0,\ \frac{d}{dt}\rho(x_i(t),x_{i+1}(t))>0
 $$
 imply that
 we have $g'(t)>0$.
 Thus at any point $A\in U^0(\Delta^{(0)})$
 we have found the vector $\vec V$ (the tangent vector to the curve
 $(x_1,x_2(\alpha_2t),\dots,x_p(\alpha_pt))$)
 such that the derivative of the function
 $g$ along $\vec V$ is greater than $0$.
 Consequently $dg(A)\ne 0$.
 Suppose we decrease $\varepsilon$.
 Note that we do not need to decrease the constructed neighborhood
 $U^0(\Delta^{(0)})$ since the critical points of $g$
 cannot appear inside it.

 Suppose we find $\varepsilon$ for a neighborhood
 $U^{p-l-1}(\Delta^{(p-l-1)})$.
 Replace this $\varepsilon$ by $\frac{\varepsilon}{2}$.
 We need to consider not the entire
 $\Delta^{(p-l)}$, but only
 $\Delta^{(p-l)}-U^{p-l-1}(\Delta^{(p-l-1)})$, i.~e.,
 we can assume that any point $A\in\Delta^{(p-l)}$
 is of the form
 $(x_1,\dots,x_1,x_{l+1},\dots,x_p)$, where
 $\rho(x_1,x_{l+1})>\varepsilon,\dots,\rho(x_1,x_p)>\varepsilon$.
 Arguing as above, we see that the whole neighborhood
 $\Delta^{p-l}$ can be obtained if
 we fix $x_1,x_{l+1},\dots,x_p$ and put $x_2,\dots,x_l$
 on geodesics passing through $x_1$.
 Again we have $g(x_1,\dots,x_p)=g(t)$, and
 $$
 g'(t)=
 \underbrace{\left(\prod_{i=l+1}^{p}\varphi(\rho(x_i,x_{i+1}))\right)}_{const}
 \left(
 \left(\prod_{i=1}^{l}\varphi(\rho(x_i,x_{i+1}))\right)'f(t)+
 \left(\prod_{i=1}^{l}\varphi(\rho(x_i,x_{i+1}))\right)f'(t)
 \right).
 $$
 Note that
 $$
 f'(t)=\underbrace{\frac{d}{dt}\sum_{i=1}^{p-1}\rho(x_i,x_{i+1})}_{\ge\alpha_l}
 +\underbrace{\frac{d}{dt}\rho(x_{l}(\alpha_lt),x_{l+1})}_{\ge-\alpha_l}\ge
 0.
 $$
 Indeed, the distance from a fixed point to a point moving
 along a geodesic cannot change with a velocity greater than $1$.
 Thus $g'(t)>0$ again.

 So we can pass from $\Delta^{(p-l-1)}$ to $\Delta^{(p-l)}$.
 Since $\Delta^{(p-2)}=\Delta$, this completes the proof.
 \proofend

 \section{Periodic billiard trajectories in a circle}

 \begin{lemma} There exists an embedding of the circle $S^1$
 into the plane $\R^2$ such that the function
 $$
 f(x_1,\dots,x_p)=\sum_{i\in\Z_p}\rho(x_i,x_{i+1})
 $$
 has $p-1$ critical points: $\frac{p-1}{2}$ maxima
 and $\frac{p-1}{2}$ points of Morse index $p-1$.
 \end{lemma}

 \begin{figure}[ht]
 \unitlength=10pt
 \begin{center}
 \begin{picture}(10,10)(0,-10)
  \drawline[-50](5,0)(8,-10)
  \drawline[-50](5,0)(2,-10)
  \drawline[-50](10,-4)(2,-10)
  \drawline[-50](0,-4)(8,-10)
  \drawline[-50](0,-4)(10,-4)
  \closecurve(5,0, 10,-4, 8,-10, 2,-10, 0,-4)
  \put(5,0){\line(5,-4){5}}
  \put(10,-4){\line(-1,-3){2}}
  \put(8,-10){\line(-1,0){6}}
  \put(2,-10){\line(-1,3){2}}
  \put(0,-4){\line(5,4){5}}
 \end{picture}\\[5mm]
 Two of the four periodic billiard trajectories for $p=5$
 \end{center}
 \end{figure}

 \beginproof
 This embedding in polar coordinates is given by the formula
 $$
 r=1-\varepsilon\cos p\varphi
 $$
 for $\varepsilon$ small enough.
 Let us show this if $p=3$ (the proof for other values of $p$
 is similar).
 $3$-periodic billiard trajectories of the non-deformed
 circle $r=1$ are inscribed regular triangles.
 The coordinates of the vertices of such triangle are
 $$
 \varphi_1=\alpha_0,\ \varphi_2=\alpha_0+\frac{2\pi}{3},
 \ \varphi_3=\alpha_0+\frac{4\pi}{3}.
 $$
 Thus $3$-periodic billiard trajectories of
 the deformed circle $r=1-\varepsilon\cos 3\varphi$ are
 $$
 \varphi_1=\alpha_0+\beta_1,\ \varphi_2=\alpha_0+\beta_2+\frac{2\pi}{3},
 \ \varphi_3=\alpha_0+\beta_3+\frac{4\pi}{3},
 $$
 where $\beta_1, \beta_2, \beta_3\to 0$ as $\varepsilon\to 0$.
 The length function is
 $$
 f(\varphi_1,\varphi_2,\varphi_3)=
 $$
 $$
 =\sqrt{(1-\varepsilon\cos 3\varphi_1)^2+(1-\varepsilon\cos 3\varphi_2)^2
 -2(1-\varepsilon\cos 3\varphi_1)(1-\varepsilon\cos 3\varphi_2)
 \cos (\varphi_1-\varphi_2)}+
 $$
 $$
 +\sqrt{(1-\varepsilon\cos 3\varphi_1)^2+(1-\varepsilon\cos 3\varphi_3)^2
 -2(1-\varepsilon\cos 3\varphi_1)(1-\varepsilon\cos 3\varphi_3)
 \cos (\varphi_1-\varphi_3)}+
 $$
 $$
 +\sqrt{(1-\varepsilon\cos 3\varphi_2)^2+(1-\varepsilon\cos 3\varphi_3)^2
 -2(1-\varepsilon\cos 3\varphi_2)(1-\varepsilon\cos 3\varphi_3)
 \cos (\varphi_2-\varphi_3)}.
 $$
 Its derivatives have the following form:
 $$
 \frac{\partial}{\partial\varphi_1}f(\varphi_1,\varphi_2,\varphi_3)=
 $$
 $$
 =\frac
 {3\varepsilon(1-\varepsilon\cos 3\varphi_1)\sin 3\varphi_1-
 (1-\varepsilon\cos 3\varphi_2)
 (3\varepsilon\sin 3\varphi_1\cos(\varphi_1-\varphi_2)-
 (1-\varepsilon\cos 3\varphi_1)\sin(\varphi_1-\varphi_2))}
 {\sqrt{(1-\varepsilon\cos 3\varphi_1)^2+(1-\varepsilon\cos 3\varphi_2)^2
 -2(1-\varepsilon\cos 3\varphi_1)(1-\varepsilon\cos 3\varphi_2)
 \cos (\varphi_1-\varphi_2)}}+
 $$
 $$
 +\frac
 {3\varepsilon(1-\varepsilon\cos 3\varphi_1)\sin 3\varphi_1-
 (1-\varepsilon\cos 3\varphi_3)
 (3\varepsilon\sin 3\varphi_1\cos(\varphi_1-\varphi_3)-
 (1-\varepsilon\cos 3\varphi_1)\sin(\varphi_1-\varphi_3))}
 {\sqrt{(1-\varepsilon\cos 3\varphi_1)^2+(1-\varepsilon\cos 3\varphi_3)^2
 -2(1-\varepsilon\cos 3\varphi_1)(1-\varepsilon\cos 3\varphi_3)
 \cos (\varphi_1-\varphi_3)}}.
 $$
 Substituting $\varphi_1$, $\varphi_2$, and $\varphi_3$ for their
 values in this formula and the same formulas for
 $\frac{\partial}{\partial\varphi_2}f(\varphi_1,\varphi_2,\varphi_3)$ and
 $\frac{\partial}{\partial\varphi_3}f(\varphi_1,\varphi_2,\varphi_3)$,
 we obtain
 $$
 \varphi_1=\alpha_0+\beta_1,\ \varphi_2=\alpha_0+\beta_2+\frac{2\pi}{3},
 \ \varphi_3=\alpha_0+\beta_3+\frac{4\pi}{3},
 $$
 $$
 \cos 3\varphi_1=\cos 3\alpha_0-3\beta_1\sin 3\alpha_0+\dots,
 $$
 $$
 \cos 3\varphi_2=\cos 3\alpha_0-3\beta_2\sin 3\alpha_0+\dots,
 $$
 $$
 \cos 3\varphi_3=\cos 3\alpha_0-3\beta_3\sin 3\alpha_0+\dots,
 $$
 $$
 \sin 3\varphi_1=\sin 3\alpha_0+3\beta_1\cos 3\alpha_0+\dots,
 $$
 $$
 \sin 3\varphi_2=\sin 3\alpha_0+3\beta_2\cos 3\alpha_0+\dots,
 $$
 $$
 \sin 3\varphi_3=\sin 3\alpha_0+3\beta_3\cos 3\alpha_0+\dots,
 $$
 $$
 \cos(\varphi_1-\varphi_2)=-\frac{1}{2}+\frac{\sqrt 3}{2}(\beta_1-\beta_2)
 +\dots,
 $$
 $$
 \cos(\varphi_2-\varphi_3)=-\frac{1}{2}+\frac{\sqrt 3}{2}(\beta_2-\beta_3)
 +\dots,
 $$
 $$
 \cos(\varphi_1-\varphi_3)=-\frac{1}{2}-\frac{\sqrt 3}{2}(\beta_1-\beta_3)
 +\dots,
 $$
 $$
 \sin(\varphi_1-\varphi_2)=-\frac{\sqrt 3}{2}-\frac{1}{2}(\beta_1-\beta_2)
 +\dots,
 $$
 $$
 \sin(\varphi_2-\varphi_3)=-\frac{\sqrt 3}{2}-\frac{1}{2}(\beta_2-\beta_3)
 +\dots,
 $$
 $$
 \sin(\varphi_1-\varphi_3)=\frac{\sqrt 3}{2}-\frac{1}{2}(\beta_1-\beta_3)
 +\dots
 $$
 Now we write that the derivatives of $f$ vanish:
 $$
 -2\sqrt 3\sin 3\alpha_0+\varepsilon\sqrt 3\sin 3\alpha_0+
 \beta_1(-6\sqrt 3\cos 3\alpha_0-\frac{\sqrt 3}{2}-\sin 3\alpha_0)+
 $$
 $$
 +\beta_2(\frac{\sqrt 3}{4}+\frac{1}{2}\sin 3\alpha_0)
 +\beta_3(\frac{\sqrt 3}{4}+\frac{1}{2}\sin 3\alpha_0)
 +\dots=0,
 $$
 $$
 -2\sqrt 3\sin 3\alpha_0+\varepsilon\sqrt 3\sin 3\alpha_0+
 \beta_1(\frac{\sqrt 3}{4}+\frac{1}{2}\sin 3\alpha_0)+
 $$
 $$
 +\beta_2(-6\sqrt 3\cos 3\alpha_0-\frac{\sqrt 3}{2}-\sin 3\alpha_0)+
 \beta_3(\frac{\sqrt 3}{4}+\frac{1}{2}\sin 3\alpha_0)
 +\dots=0,
 $$
 $$
 -2\sqrt 3\sin 3\alpha_0+\varepsilon\sqrt 3\sin 3\alpha_0+
 \beta_1(\frac{\sqrt 3}{4}+\frac{1}{2}\sin 3\alpha_0)+
 $$
 $$
 +\beta_2(\frac{\sqrt 3}{4}+\frac{1}{2}\sin 3\alpha_0)+
 \beta_3(-6\sqrt 3\cos 3\alpha_0-\frac{\sqrt 3}{2}-\sin 3\alpha_0)
 +\dots=0.
 $$
 Note that the constant term must be equal to $0$.
 Consequently ${\alpha_0=0}$ or ${\alpha_0=\frac{\pi}{3}}$.
 It is clear that
 $\beta_1=\beta_2=\beta_3=0$ is a periodic billiard
 trajectory for any $\varepsilon$, hence
 all coefficients of $\varepsilon^k$ are equal to $0$.
 Thus the dominant terms in this system are
 $\beta_1$, $\beta_2$, and $\beta_3$ with their coefficients.
 The linear system for
 $\beta_1$, $\beta_2$, $\beta_3$ has only the trivial solution.
 Hence we have the two trajectories:
 $(\frac{\pi}{3},\pi,\frac{5\pi}{3})$ for the maximum and
 $(0,\frac{2\pi}{3},\frac{4\pi}{3})$ for the point of index $2$.
 \proofend

 \begin{corollary} The Euler characteristic of the space
 $T^p/D_p-\Delta$ is equal to $0$.
 Moreover, we can estimate the dimensions of homology:
 $$
 \dim H_p(T^p/D_p,\Delta;\Z_2)=\dim H_{p-1}(T^p/D_p,\Delta;\Z_2)\le
 \frac{p-1}{2},
 $$
 $$
 \dim H_q(T^p/D_p,\Delta;\Z_2)=0,\ q\le p-2.
 $$
 \end{corollary}

 \begin{lemma}
 $$
 \dim H_p(T^p/D_p,\Delta;\Z_2)\ge\frac{p-1}{2}.
 $$
 \end{lemma}

 \beginproof We construct a cellular division of the torus
 $T^p$.
 The torus $T^p$ is a Cartesian product of $p$ circles.
 Let $x_i\in[0,1]$ be the cyclic coordinate on the $i$th circle.
 To each permutation
 $\sigma\in S_p$ we assign a $p$-dimensional cell
 $$
 e^p_{\sigma}=\{x_{\sigma(1)}<\dots<x_{\sigma(p)}\}.
 $$
 This division is invariant with respect to the action of the permutation
 group $S_p$.
 Then the sums of cells belonging to one connected component
 of \ $T^p/D_p-\Delta$ form a basis in
 $H_p(T^p/D_p,\Delta;\Z_2)$.
 We construct a function $I$ of a cell such that $I$
 is constant on the connected components of
 $T^p/D_p-\Delta$.

 Suppose we have a cell $e^p_{\sigma}\subset T^p$.
 For each pair of indices $(i,i+1)$, $i\in\Z_p$,
 one of the two inequalities
 $x_i<x_{i+1}$ or $x_i>x_{i+1}$ holds. By definition, put
 $$
 I(e^p_{\sigma})=|\#(>)-\#(<)|.
 $$
 Cells of the same connected component can be obtained
 from $e^p_{\sigma}$ by a sequence of the following transformations:
 $$
 (1)\ e^p_{\sigma} \mapsto e^p_{t\sigma}, t\in D_p,
 $$
 $$
 (2)\ e^p_{\sigma} \mapsto
 \{x_{\sigma(2)}<\dots<x_{\sigma(p)}<x_{\sigma(1)}\},
 $$
 $$
 (3)\ e^p_{\sigma} \mapsto e^p_{\tau_{i}\sigma},
 \tau_{i}=
 \left({{1,\dots,}\atop{1,\dots,}}
 {{i,}\atop{i+1,}}
 {{i+1,}\atop{i,}}
 {{\dots,p}\atop{\dots,p}}\right),
 \ p-2\ge|\sigma(i+1)-\sigma(i)|\ge 2,
 \ i=1,\dots,p-1.
 $$
 Evidently, transformations
 (1) and (3) do not change
 $I(e^p_{\sigma})$. Transformation (2) does not change it either:
 $$
 {x_{\sigma(1)}<x_{\sigma(1)+1}\atop
 x_{\sigma(1)-1}>x_{\sigma(1)}}\longrightarrow
 {x_{\sigma(1)}>x_{\sigma(1)+1}\atop
 x_{\sigma(1)-1}<x_{\sigma(1)}}
 $$
 Now we find $\frac{p-1}{2}$ $p$-cells with different
 values of $I$:
 $$
 \begin{array}{c}
 \{x_1<\dots<x_p\}\\
 \{x_1<x_3<x_2<x_4<\dots<x_p\}\\
 \{x_1<x_4<x_3<x_2<x_5<\dots<x_p\}\\
 \vdots\\
 \{x_1<x_{\frac{p+1}{2}}<\dots<x_2<x_{\frac{p+1}{2}+1}\ \dots\ <x_p\}
 \end{array}
 $$
 The function $I$ takes values $p-2,p-4,\dots,1$ at these cells.
 Thus $T^p/D_p-\Delta$ has at least
 $\frac{p-1}{2}$ connected components.
 The lemma is proved.
 \proofend

 \medskip
 So we obtain the following result:

 \begin{theorem}
 $$
 BT_p(S^1)=p-1.
 $$
 \end{theorem}

 \begin{note} George Birkhoff found this estimate,
 but he did not prove that it is exact.
 \end{note}

 \section{Periodic billiard trajectories in almost round spheres}

 Consider the standard $n$-sphere embedded in $\R^{n+1}$:
 $$
 S^n=\{(y_0,y_1,\dots,y_n)\in\R^{n+1}: \sum{y_i^2}=1\}
 $$
 The length function $f(x_1,\dots,x_p)$ has
 $\frac{p-1}{2}$ diffeomorphic non-degenerate critical manifolds.
 Each critical manifold is the set of regular
 $n$-gons (or "stars") inscribed in the unit circle
 centered at $0\in\R^{n+1}$.
 Denote such a component by $V_{n,p}$.

 \begin{lemma}
 $$
 \sum_{q=0}^{\infty}\dim H_q(V_{n,p};\Z_2)=2n
 $$
 \end{lemma}

 \beginproof It is clear that $V_{n,p}$
 is a bundle over the Grassmannian $G_{2,n+1}$
 with fiber $S^1$.
 The fiber has one-dimensional homologies,
 so this bundle is $\Z_2$-homology simple.
 Consider the
 $\Z_2$-cohomology spectral sequence of this bundle.
 Differentials of $E_2$ are multiplications by
 the characteristic class
 $\alpha\in H^2(G_{2,n+1};\Z_2)$ of this bundle.
 Grassmannian $G_{2,n+1}$ is subdivided into Schubert cells
 $\sigma_{2,\dots,2,1,\dots,1}$ --- see~\cite{griffiths}.
 In this cell complex all boundary operators vanish,
 therefore each cell corresponds to a class of homology.
 By the same symbol
 $\sigma_{2,\dots,2,1,\dots,1}$ denote the Poincar\'{e} dual
 cohomology class.
 Note that
 $$
 \dim H^q(G_{2,n+1};\Z_2)=\left\{ {[\frac q2]+1, \mbox{ if }q\le n-1
 \atop
  [\frac{2(n-1)-q}{2}]+1, \mbox{ if }q>n-1} \right.
 $$
 Indeed, a number $q$ can be presented as a sum
 of ones and twos in
 $[\frac q2]+1$ ways.
 If $q\le n-1$, then each sum corresponds to
 a Schubert cell of the Grassmannian $G_{2,n+1}$.
 If $q>n-1$, then our assertion follows from Poincar\'{e} duality.
 Moreover, we have the Pieri formula:
 $$
 \sigma_a\smallsmile\sigma_{b_1,b_2,\dots}=
 \sum_{b_i\le c_i\le b_{i-1}\atop \sum c_i=a+\sum b_i}\sigma_{c_1,c_2,\dots}
 $$

 Let us show that the characteristic class
 of our bundle is equal to $\sigma_2$.
 First suppose $n=2$, $p=3$.
 Then $G_{2,n+1}$ is the
 projective plane $\RP^2$.
 $V_{2,3}$ is the manifold of all big regular
 triangles inscribed in the unit sphere.
 Points of the base $\RP^2$ are the lines orthogonal
 to the planes of these triangles.
 $\RP^2$ is divided into three cells
 $e^2$, $e^1$, and $e^0$ of dimensions $2$, $1$, and $0$.
 The cell $e^1$ consists of all lines
 that belong to the coordinate plane
 $Oxy$ (except the axis $Ox$).
 All other lines belong to the cell $e_2$.
 Each of them is defined by a point of the upper hemisphere.
 To compute the characteristic class we
 must construct a section of the bundle
 $s:e^1\to V_{2,3}$.
 Suppose this section consists of all vertical triangles,
 i.e., triangles with one vertex at the north pole of the sphere.
 By $h:\overline{D^2}\to\RP^2$ denote
 the characteristic map of the cell $e^2$.
 We assume that the open disk
 $D^2$ is the upper hemisphere.
 We have coordinates
 $(x,y)$ on the disk.
 The bundle over the disk is trivial.
 Introduce coordinates of direct product in it.
 Suppose $\delta$ is a triangle of $V_{2,3}$
 not lying in any vertical plane.
 The plane of this triangle and the plane
 $Oxy$ intersect by the line $l$.
 Rotate this triangle around the line $l$ to the
 horizontal position (see the picture).

 \begin{figure}[ht]
 \unitlength=1cm
 \begin{center}
 \begin{picture}(10,10)(0,0)
 \put(0,0){\epsfxsize=10cm\epsfbox{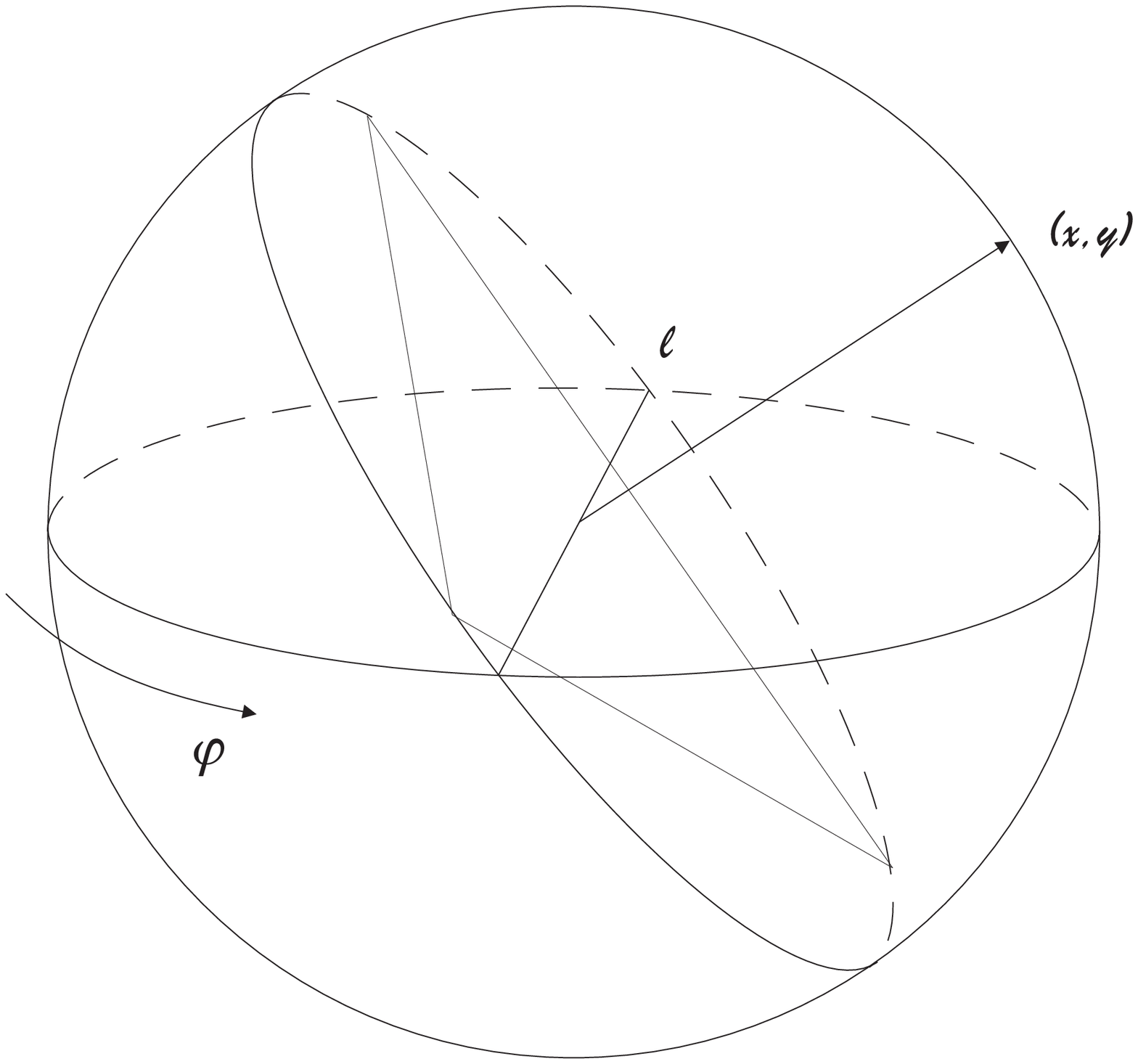}}
 \end{picture}
 \end{center}
 \end{figure}

 Let $\varphi\in[0,\frac{2\pi}{3}]$ be
 the smallest polar angle of its vertices.
 Then the coordinates
 $(x,y,\varphi)$ trivialize the bundle over the disk $D^2$.
 Now consider a bundle induced by the
 characteristic map $h$ over the
 closed disk $\overline{D^2}$.
 We can assume that any element of this fibered space
 is the big triangle with a normal that looks to
 the upper hemisphere or lies in the plane $Oxy$.
 This space has the same coordinates $(x,y,\varphi)$.
 We have a section
 $s':S^1_{\varphi'}\to\overline{D^2}\times S^1_{\varphi}$
 over the boundary of the disk,
 where $\varphi=\varphi'\bmod\frac{2\pi}{3}$ ---
 if we go around the circle--boundary once,
 then we make three complete turns of the circle--fiber.
 Thus the $\Z_2$-characteristic class of this
 bundle assigns the number $1$ to the cell $e^2$.
 Similarly for $p>3$ this class is also equal to
 $\sigma_2$.
 If $n>2$, then we have an embedding $G_{2,3}\to G_{2,n+1}$
 that preserves the cellular division.
 There is one $2$-cell in $G_{2,3}$ and two $2$-cells
 in $G_{2,n+1}$ for $n>2$:
 $\sigma_{2,\dots,2}$ (the image of the $2$-cell of $G_{2,3}$)
 and $\sigma_{2,\dots,2,1,1}$.
 It is clear that the cell $\sigma_{2,\dots,2,1,1}$
 does not make any obstruction for the section over
 $1$-cell $\sigma_{2,\dots,2,1}$.
 Thus in this case the characteristic class is also equal
 to $\sigma_2$.

 We have in $H^*(G_{2,n+1};\Z_2)$:
 $$
 \sigma_2\smallsmile\sigma_{b_1,\dots,b_k}=\sigma_{2,b_1,\dots,b_k},
 \mbox{ if } k<n-1,
 $$
 $$
 \sigma_2\smallsmile\sigma_{b_1,\dots,b_{n-1}}=0.
 $$
 by the Pieri formula.
 Now we see that all differentials
 $d^{j,1}_{2}:E^{j,1}_{2}\to E^{j+2,0}_2$
 are either monomorphisms or epimorphisms,
 $d^{n-2,1}_2$ is an isomorphism.
 Thus the term
 $E_3=E_{\infty}$ is
 $$
 \begin{array}{ccccccccc}
 0 & 0 & \ldots & 0 & \Z_2 & \Z_2 & \ldots & \Z_2 & \Z_2 \\
 \Z_2 & \Z_2 & \ldots & \Z_2 & \Z_2 & 0 & \ldots & 0 & 0 \\
 \end{array}
 $$
 All $H^q(V_{n,p};\Z_2)=\Z_2$ for $q=0,\dots,2n-1$.
 This completes the proof.
 \proofend

 \begin{note}
 Suppose $f$ is a function with a non-degenerate
 critical manifold $V\subset\{f=a\}$.
 Then its small Morse deformation
 $\tilde{f}$ has at least
 $\sum\dim H_q(V;\Z_2)$ critical points.
 \end{note}
 \beginproof
 The function $\tilde{f}\mid_V$ is Morse in the general case.
 Suppose $\tilde{f}\mid_V$ has $N$ critical points.
 Then $N\ge\sum\dim H_q(V;\Z_2)$.
 A point $A\in V$ is critical for the function $\tilde{f}\mid_V$
 if $\grad\tilde{f}\perp T_AV$.
 In a neighborhood of a critical point we have
 $$
 f=x_1^2+\dots+x_{\alpha}^2-x_{\alpha+1}^2-x_k^2,
 $$
 $$
 \tilde{f}=x_1^2+\dots+x_{\alpha}^2-x_{\alpha+1}^2-x_k^2+
 \sum_{i=1}^{m}\varepsilon_ix_i+o(\mbox{\bf x}),
 $$
 $$
 V=\{x_{k+1}=\dots=x_m=0\}.
 $$
 Since $\grad\tilde{f}\perp V$, we have that
 $\varepsilon_{k+1}=\dots=\varepsilon_m=0$.
 Hence the function
 $\tilde{f}$ has only one critical point in this neighborhood:
 $x_i\approx\frac{\varepsilon_i}{2}$.
 Thus the function
 $\tilde{f}$ has at least $N$ critical points.
 \proofend

 We have proved the following fact.
 \begin{theorem}
 Suppose $p$ is an odd prime, $n$ is any integer.
 Then a generic small perturbation of a standard
 round $n$-sphere has at least $n(p-1)$
 $p$-periodic billiard trajectories.
 \end{theorem}

 \section{$3$-periodic billiard trajectories in a $2$-dimensional sphere}

 \begin{theorem}
 The minimal number of $3$-periodic billiard trajectories
 in the sphere satisfies
 $$
 BT_3(S^2)\ge 4.
 $$
 \end{theorem}
 \beginproof
 We shall construct a cellular division of the space
 $S^2\times S^2\times S^2/D_3$.
 We assume that our sphere is a square
 with the boundary contracted to a point:
 $$
 S^2=[0,1]_{\varphi}\times[0,1]_{\psi}/\{\varphi=0,\varphi=1,\psi=0,\psi=1\}
 $$
 The list of all cells is
 $$
 \dim=6
 $$
 $$
 \omega_1=\left\{ {\varphi_1<\varphi_2<\varphi_3 \atop
                   \psi_1<\psi_2<\psi_3} \right\}\
 \omega_2=\left\{ {\varphi_1<\varphi_2<\varphi_3 \atop
                   \psi_1<\psi_3<\psi_2} \right\}\
 \omega_3=\left\{ {\varphi_1<\varphi_2<\varphi_3 \atop
                   \psi_2<\psi_1<\psi_3} \right\}\
 $$
 $$
 \omega_4=\left\{ {\varphi_1<\varphi_2<\varphi_3 \atop
                   \psi_2<\psi_3<\psi_1} \right\}\
 \omega_5=\left\{ {\varphi_1<\varphi_2<\varphi_3 \atop
                   \psi_3<\psi_1<\psi_2} \right\}\
 \omega_6=\left\{ {\varphi_1<\varphi_2<\varphi_3 \atop
                   \psi_3<\psi_2<\psi_1} \right\}\
 $$
 $$
 \dim=5
 $$
 $$
 \sigma_1=\left\{ {\varphi_1<\varphi_2<\varphi_3 \atop
                   \psi_1=\psi_2<\psi_3} \right\}\
 \sigma_2=\left\{ {\varphi_1<\varphi_2<\varphi_3 \atop
                   \psi_3<\psi_1=\psi_2} \right\}\
 \sigma_3=\left\{ {\varphi_1<\varphi_2<\varphi_3 \atop
                   \psi_1=\psi_3<\psi_2} \right\}\
 $$
 $$
 \sigma_4=\left\{ {\varphi_1<\varphi_2<\varphi_3 \atop
                   \psi_2<\psi_1=\psi_3} \right\}\
 \sigma_5=\left\{ {\varphi_1<\varphi_2<\varphi_3 \atop
                   \psi_2=\psi_3<\psi_1} \right\}\
 \sigma_6=\left\{ {\varphi_1<\varphi_2<\varphi_3 \atop
                   \psi_1<\psi_2=\psi_3} \right\}\
 $$
 $$
 \delta_1=\left\{ {\varphi_1=\varphi_2<\varphi_3 \atop
                   \psi_1<\psi_2<\psi_3} \right\}\
 \delta_2=\left\{ {\varphi_3<\varphi_1=\varphi_2 \atop
                   \psi_1<\psi_2<\psi_3} \right\}\
 \delta_3=\left\{ {\varphi_1=\varphi_3<\varphi_2 \atop
                   \psi_1<\psi_2<\psi_3} \right\}\
 $$
 $$
 \delta_4=\left\{ {\varphi_2<\varphi_1=\varphi_3 \atop
                   \psi_1<\psi_2<\psi_3} \right\}\
 \delta_5=\left\{ {\varphi_2=\varphi_3<\varphi_1 \atop
                   \psi_1<\psi_2<\psi_3} \right\}\
 \delta_6=\left\{ {\varphi_1<\varphi_2=\varphi_3 \atop
                   \psi_1<\psi_2<\psi_3} \right\}\
 $$
 $$
 \dim=4
 $$
 $$
 \alpha_1=\left\{ {\varphi_1=\varphi_2=\varphi_3 \atop
                   \psi_1<\psi_2<\psi_3} \right\}\
 \alpha_2=\left\{ {\varphi_1<\varphi_2<\varphi_3 \atop
                   \psi_1=\psi_2=\psi_3} \right\}\
 \beta_1=\left\{ {\varphi_1=\varphi_2<\varphi_3 \atop
                   \psi_1<\psi_2=\psi_3} \right\}\
 $$
 $$
 \beta_2=\left\{ {\varphi_3<\varphi_1=\varphi_2 \atop
                   \psi_1<\psi_2=\psi_3} \right\}\
 \beta_3=\left\{ {\varphi_1=\varphi_2<\varphi_3 \atop
                   \psi_2=\psi_3<\psi_1} \right\}\
 \beta_4=\left\{ {\varphi_3<\varphi_1=\varphi_2 \atop
                   \psi_2=\psi_3<\psi_1} \right\}\
 $$
 $$
 \gamma_1=\left\{ {\varphi_1=0,\varphi_2<\varphi_3 \atop
                   \psi_1=0,\psi_2<\psi_3} \right\}\
 \gamma_2=\left\{ {\varphi_1=0,\varphi_2<\varphi_3 \atop
                   \psi_1=0,\psi_3<\psi_2} \right\}\
 $$
 $$
 \dim=3
 $$
 $$
 \kappa_1=\left\{ {\varphi_1=0,\varphi_2=\varphi_3 \atop
                   \psi_1=0,\psi_2<\psi_3} \right\}\
 \kappa_2=\left\{ {\varphi_1=0,\varphi_2<\varphi_3 \atop
                   \psi_1=0,\psi_2=\psi_3} \right\}\
 $$
 Now we compute the boundary operators.
 $$
 \partial_6:C_6\to C_5
 $$
 $$
 \omega_1\mapsto\delta_1+\delta_6+\sigma_1+\sigma_6,\ \
 \omega_2\mapsto\delta_3+\delta_6+\sigma_3+\sigma_6,\ \
 $$
 $$
 \omega_3\mapsto\delta_1+\delta_4+\sigma_1+\sigma_4,\ \
 \omega_4\mapsto\delta_2+\delta_3+\sigma_4+\sigma_5,\ \
 $$
 $$
 \omega_5\mapsto\delta_4+\delta_5+\sigma_2+\sigma_3,\  \
 \omega_6\mapsto\delta_2+\delta_5+\sigma_2+\sigma_5\  \
 $$
 $$
 \partial_5:C_5\to C_4
 $$
 $$
 \sigma_1\mapsto\alpha_2+\beta_4,\ \
 \sigma_2\mapsto\alpha_2+\beta_2,\ \
 \sigma_3\mapsto\alpha_2+\beta_3+\beta_4,\ \
 $$
 $$
 \sigma_4\mapsto\alpha_2+\beta_1+\beta_2,\ \
 \sigma_5\mapsto\alpha_2+\beta_3,\ \
 \sigma_6\mapsto\alpha_2+\beta_1\ \
 $$
 $$
 \delta_1\mapsto\alpha_1+\beta_1,\ \
 \delta_2\mapsto\alpha_1+\beta_2,\ \
 \delta_3\mapsto\alpha_1+\beta_1+\beta_3,\ \
 $$
 $$
 \delta_4\mapsto\alpha_1+\beta_2+\beta_4,\ \
 \delta_5\mapsto\alpha_1+\beta_3,\ \
 \delta_6\mapsto\alpha_1+\beta_4,\ \
 $$
 $$
 \partial_4:C_4\to C_3
 $$
 $$
 \alpha_1, \alpha_2, \beta_1, \beta_2, \beta_3, \beta_4\mapsto 0,\ \
 \gamma_1, \gamma_2\mapsto\kappa_1+\kappa_2
 $$
 $$
 \partial_3=0
 $$
 Now we compute the kernels and the images
 of the boundary operators and the homology groups.
 $$
 \ker\partial_6=\langle\omega_1+\omega_2+\omega_3+\omega_4+
 \omega_5+\omega_6\rangle,\ \
 \im\partial_7=0,\ \ H_6=\Z_2,
 $$
 $$
 \ker\partial_5=\langle\sigma_1+\sigma_2+\sigma_3+
 \sigma_4+\sigma_5+\sigma_6\rangle
 \oplus\im\partial_6,\ \ H_5=\Z_2,
 $$
 $$
 \im\partial_5=\langle\alpha_1,\alpha_2,\beta_1,\beta_2,\beta_3,\beta_4\rangle,\
 \
 \ker\partial_4=\langle\gamma_1+\gamma_2\rangle\oplus\im\partial_5,\ \
 H_4=\Z_2,
 $$
 $$
 \ker\partial_3=\langle\gamma_1,\gamma_2\rangle,\ \
 \im\partial_4=\langle\gamma_1+\gamma_2\rangle,\ \
 H_3=\Z_2.
 $$
 Thus we obtain that
 $$
 H_*(S^2\times S^2\times S^2/D_3,\Delta;\Z_2)=
 \{0,0,0,\Z_2,\Z_2,\Z_2,\Z_2\}.
 $$
 This completes the proof.
 \proofend

 \section{A general estimate for $3$-periodic billiard trajectories}

 \begin{lemma}
 Suppose we have an exact sequence of vector spaces
 $$
 0\stackrel{i_m}{\longrightarrow}B_m
 \stackrel{j_m}{\longrightarrow}C_m
 \stackrel{\partial_m}{\longrightarrow}A_{m-1}
 \stackrel{i_{m-1}}{\longrightarrow}B_{m-1}
 \stackrel{j_{m-1}}{\longrightarrow}C_{m-1}
 \stackrel{\partial_{m-1}}{\longrightarrow}\ldots
 $$
 By definition, put $a_q=dim A_q,\ b_q=\dim B_q,\ c_q=\dim C_q.$
 Then $\sum c_q\ge\sum b_q-\sum a_q$.
 \end{lemma}

 \beginproof
 $$
 b_m=\dim\im j_m=
 \dim\ker\partial_m=
 c_m-\dim\im\partial_m=
 c_m-\dim\ker i_{m-1}=
 $$
 $$
 =c_m-(a_{m-1}-\dim\im i_{m-1})\le
 c_m+(a_{m-1}-\dim\im i_{m-1})=
 $$
 $$
 =c_m+a_{m-1}-\dim\ker j_{m-1}=
 c_m+a_{m-1}-b_{m-1}+\dim\im j_{m-1}=
 $$
 $$
 =c_m+a_{m-1}-b_{m-1}+\dim\ker\partial_{m-1}=
 $$
 $$
 =c_m+a_{m-1}-b_{m-1}+c_{m-1}-\dim\im\partial_{m-1}\le
 $$
 $$
 \le c_m+a_{m-1}-b_{m-1}+c_{m-1}+a_{m-2}-\dim\im i_{m-2}\le\ldots\le
 $$
 $$
 \le \sum_{q=0}^{m}c_q+\sum_{q=0}^{m-1}a_q-\sum_{q=0}^{m-1}b_q\ .
 $$
 \proofend

 \begin{theorem}
 Suppose $M$ is a closed manifold,
 $B=\sum\dim H_q(M;\Z_3)$.
 Then we have an estimate
 $$
 BT_3(M)\ge\frac{B^3-3B^2+2B}{6}.
 $$
 \end{theorem}
 \beginproof
 The group $\Z_3$ acts on $M\times M\times M$:
 $$
 t:M\times M\times M\to M\times M\times M,\ \ t(x_1,x_2,x_3)=(x_2,x_3,x_1),
 \ \ t^3=id.
 $$
 $\Delta^{(0)}=\{(x,x,x)\}\subset M\times M\times M$ is the
 set of all fixed points of this action.
 $X=M\times M\times M/\Z_3$.
 The Smith theory gives us an estimate
 (see \cite{viro}):
 $$
 \dim H_q(X,\Delta^{(0)};\Z_3)\ge\frac{\dim H_q(M\times M\times M;\Z_3)-
 \dim H_q(\Delta^{(0)};\Z_3)}{3}.
 $$
 As above, $\Delta=\{(x,x,y)\}\subset X$ is the diagonal
 of the cyclic cube $M\times M\times M/\Z_3$.
 Now we must estimate the dimensions of
 the other relative homology groups
 $H_q(X,\Delta)$.
 First we have an exact homological sequence of the pair
 $(\Delta,\Delta^{(0)})$.
 The pair $(\Delta,\Delta^{(0)})$ coincides with the pair
 $(M\times M, M=\{(x,x)\}\subset M\times M)$, and we have
 $$
 \ldots\longrightarrow H_q(M)
 \stackrel{i_*}{\longrightarrow}
 H_q(M\times M)
 \stackrel{j_*}{\longrightarrow}
 H_q(M\times M,M)
 \stackrel{\partial_*}{\longrightarrow}
 H_{q-1}(M)
 \longrightarrow\ldots
 $$
 Let us note that $\ker i_*=\{0\}$.
 Indeed, suppose $p:M\times M\to M$ is a projection
 on any of the two factors.
 Then $p\circ i=id$, hence
 $p_*\circ i_*=id_*$, i. e, $\ker i_*=\{0\}$. Further,
 $$
 H_q(\Delta)\cong H_q(\Delta,\Delta^{(0)})\oplus H_q(\Delta^{(0)}).
 $$
 So it is clear that
 $$
 \sum_q\dim H_q(\Delta,\Delta^{(0)})=B^2-B.
 $$
 Now consider the exact sequence of the triple
 $(X,\Delta,\Delta^{(0)})$:
 $$
 \ldots\longrightarrow H_q(\Delta,\Delta^{(0)})
 \longrightarrow H_q(X,\Delta^{(0)})
 \longrightarrow H_q(X,\Delta)
 \longrightarrow H_{q-1}(\Delta,\Delta^{(0)})
 \longrightarrow\ldots
 $$
 Therefore
 $$
 \sum_q\dim H_q(X,\Delta)\ge\sum_q\dim H_q(X,\Delta^{(0)})-
 \sum_q\dim H_q(\Delta,\Delta^{(0)}).
 $$
 Finally we use the modified Theorem 1:
 $$
 BT_p(M)\ge\frac{1}{2}
 \sum_{q=1}^{mp}\dim H_q(M\times\dots\times M/\Z_p,\Delta;\Z_p)
 $$
 This completes the proof.
 \proofend

 \end{document}